# QUANTILE REGRESSION WITH VARYING COEFFICIENTS

By Mi-Ok Kim

*University of Kentucky*

Quantile regression provides a framework for modeling statistical quantities of interest other than the conditional mean. The regression methodology is well developed for linear models, but less so for nonparametric models. We consider conditional quantiles with varying coefficients and propose a methodology for their estimation and assessment using polynomial splines. The proposed estimators are easy to compute via standard quantile regression algorithms and a stepwise knot selection algorithm. The proposed Rao-score-type test that assesses the model against a linear model is also easy to implement. We provide asymptotic results on the convergence of the estimators and the null distribution of the test statistic. Empirical results are also provided, including an application of the methodology to forced expiratory volume (FEV) data.

**1. Introduction.** Quantile regression has appeared as an alternative to least squares in a wide range of applications. When the center of the conditional distribution of a response variable $Y$, given a covariate vector $X$, is under investigation, median regression provides a consistent estimator of the conditional median without assuming a specific form for the conditional distribution. When other conditional quantities, for example, the lower or upper tail of the conditional distribution, are of interest, quantile regression provides a way to directly estimate the interesting quantities without assuming that such quantiles are related to $X$ in the same fashion as the conditional mean. A good example can be found in a study of consumer demand for electricity, where heavy users responded much more drastically to weather and time variation than average users [17].

Quantile regression has been implemented with nonparametric methods to overcome the limitation of a linear model. Local polynomial, smoothing spline and B-spline smoothing methods were considered in [4, 35], [13, 24] and [14], respectively. He and Shi [15] used bivariate tensor product B-splines with a partly linear model. Wei and He [34] considered B-splines with a partly linear model for longitudinal data to model growth curves. De Gooijer and Zerom [6] and Doksum and Koo [7] considered local polynomials and polynomial splines with additive models. Nevertheless, compared to its development with a linear model, the nonparametric implementation of quantile regression has been limited.

In this paper, we consider a varying-coefficient model for the conditional quantiles: for a random vector $(T, X, Y) \in [0,1] \times \mathbb{R}^{p+1} \times \mathbb{R}$, we suppose that

$$(1.1) \qquad q_\tau(t, x) = \beta_0(t)x_0 + \beta_1(t)x_1 + \cdots + \beta_p(t)x_p,$$

where $q_\tau(t, x)$ denotes the $\tau$th conditional quantile of $Y$, given that $(T, X) = (t, x)$, and the $\beta_j(t)$ are unknown smooth functions of $t$ for $j = 0, 1, \ldots, p$. A random sample is modeled by

$$(1.2) \qquad Y_i = q_\tau(T_i, X_i) + e_i, \qquad 1 \le i \le n,$$

where the $e_i$ are independent random variables with $\tau$th quantile 0 and are independent of the $(T_i, X_i)$.

Varying-coefficient models constitute an important class of nonparametric models: they have been widely applied to analyze conditional means due to their flexibility and interpretability (see [2], [9] and [19] for examples of the analysis of time series, survival and longitudinal data). However, a varying-coefficient model for the conditional quantiles was considered only recently by [18] and [3], where local polynomials were used for independent and time series data, respectively. In this paper, we provide a more comprehensive treatment of the model for independent data by proposing a polynomial-spline-based methodology for its estimation and assessment.

The proposed methodology is quite general, in the sense that we do not require the unknown coefficients $\beta_j(t)$ to be smooth curves of a common degree or the $e_i$ to be identically distributed; such assumptions facilitate the development of asymptotic theories covered in the later sections. The methodology is also simple and easy to implement: the proposed estimator is readily computed by standard quantile regression algorithms and the proposed stepwise knot selection algorithm. Different knots can be selected for different coefficients and the algorithm can accommodate varying degrees of smoothness without complicating the computation (see [11] and [19] for discussion of this issue with respect to local polynomials and smoothing splines). We show that the spline estimators and their derivatives attain the optimal rates of global convergence under appropriate conditions.

The stated feature of accommodating different degrees of smoothness in the coefficient functions is shared by any estimation method employing polynomial splines, such as the one which appears in [21]. However, the one in [21] requires its knots to be equispaced and selects the numbers of knots for the coefficients via cross-validation. In addition, it uses the least squares method for estimation since it concerns the conditional mean.

For the assessment of the model against a linear model, a Rao-score-type test is proposed. The test uses polynomial spline smoothing to account for the structure of the varying-coefficient design space and is therefore potentially more powerful than a test with a nonspecific nonparametric alternative. Implementation of the test is easy, since the test statistic has a simple asymptotic null distribution.

Model assessment is essential in this context since a linear model will be preferred on the grounds of parsimony unless the underlying science warrants use of the varying-coefficient model. While this issue is generally important for nonparametric models, it has received little attention in quantile regression. Only recently, He and Zhu [16] and Horowitz and Spokoiny [20] considered general lack-of-fit tests for linear quantile regression. Both tests are consistent for any fixed alternative. However, they are not powerful in this context because their alternative space does not take into account the structure of the varying-coefficient model. Any lack-of-fit test with a nonspecific alternative will have the same problem and is thus not considered here (see [12] for a survey of such general tests for the conditional mean).

Likelihood-ratio-type tests were proposed for the conditional mean in a similar context [1, 2, 10]. Local polynomials were used to smooth the varying-coefficient alternative design space. In essence, these tests compared the residual sums of squares under the null and the alternative. We adapted and used the tests to evaluate the performance of the proposed test in a Monte Carlo study. The results are presented in Section 4.1.

We present a motivating example involving forced expiratory volume (FEV) data collected on children aged three to nineteen [29]. A topic of interest is the effect of cigarette smoking on FEV after accounting for the effects of age, sex and height. One may speculate that the effects of the independent variables change with age. On the other hand, the lower conditional quantiles of FEV are of independent interest as they can provide a gauge for poor pulmonary function. Both research equations can be investigated by fitting model (1.2) and testing for age dependency in the effects of the independent variables. An important point is that we do not need to assume that the effects of the independent variables are same across the quantiles. Detailed results for this example will be given in Section 4.2.

In this paper, asymptotic results will be proven under rather simplified conditions that require knots for the splines to be nonstochastic and quasi-uniformly placed. As with most other asymptotic analyses in the literature, the results obtained under such conditions can provide insight into the large-sample behavior of the proposed methodology in more general settings. The rest of the article is organized as follows. Sections 2 and 3 discuss model estimation and assessment. Section 4 presents empirical results. Proofs are deferred to the Appendix.

## 2. Estimation.

2.1. *Regression splines.* The degree of smoothness of the true coefficient functions determines how well the functions can be approximated. We assume that the $\beta_j(t)$ are functions of $t$ with a common degree of smoothness $r$ defined as follows: let $\mathcal{H}_r$ be the collection of all functions on $[0,1]$ for which the $m$th order derivative satisfies the Hölder condition of order $\gamma$ with $r \equiv m + \gamma$. That is, for each $h \in \mathcal{H}_r$, $|h^{(m)}(s) - h^{(m)}(t)| \leq W_0|s-t|^\gamma$ for any $0 \leq s, t \leq 1$ and a positive finite constant $W_0$.

CONDITION 1. $\beta_j(t) \in \mathcal{H}_r$, $j = 0, 1, \ldots, p$, for some $r > 1/2$.

If the $\beta_j$ have bounded $d$th order derivatives on $[0,1]$, then Condition 1 holds with $r = d$.

We shall use normalized B-splines of order $m+1$ to approximate the $\beta_j(t)$. We consider a sequence of positive integers $\{k_n\}$, $n \geq 1$, and an extended partition of $[0,1]$ by $k_n$ quasi-uniform knots. Following [30], we denote the associated B-spline basis functions by $B_j(t)$, $j = 1, \ldots, k_n + m + 1$. We define $\pi_{k_n}(t) = (B_1(t), \ldots, B_{k_n+m+1}(t))^{\mathrm{T}}$ and $\Pi_{k_n}(t,x) = (x_0\pi_{k_n}(t)^{\mathrm{T}}, x_1\pi_{k_n}(t)^{\mathrm{T}}, \ldots, x_p\pi_{k_n}(t)^{\mathrm{T}})^{\mathrm{T}}$.

The proposed polynomial spline estimator of $q_\tau(t,x)$ is given by

$$(2.1) \qquad \hat{q}_{k_n}(t,x) = \sum_{j=0}^p \hat{\theta}_j^{\mathrm{T}} \pi_{k_n}(t) x_j = \hat{\boldsymbol{\theta}}_{k_n}^{\mathrm{T}} \Pi_{k_n}(t,x),$$

where $p_{k_n} = (p+1)(k_n + m + 1)$ and $\hat{\boldsymbol{\theta}}_{k_n} = (\hat{\theta}_0^{\mathrm{T}}, \ldots, \hat{\theta}_p^{\mathrm{T}})^{\mathrm{T}} \in \mathbb{R}^{p_{k_n}}$ solves the minimization problem

$$(2.2) \qquad \min_{\boldsymbol{\theta} \in \mathbb{R}^{p_{k_n}}} \sum_{1 \leq i \leq n} \rho_\tau(Y_i - \boldsymbol{\theta}^{\mathrm{T}} \Pi_{k_n}(T_i, X_i))$$

for $\rho_\tau(s) = s(\tau - I(s < 0))$. Accordingly, the polynomial spline estimator $\hat{\beta}_j(t)$ of $\beta_j(t)$ is given by $\hat{\theta}_j^{\mathrm{T}} \pi_{k_n}(t)$ for each $j$.

2.2. *Asymptotic result.* The following are sufficient conditions for the proposed polynomial spline estimator and its derivatives to converge at their best possible rates as the sample size goes to infinity. For simplicity, we adopt the following notation throughout the paper: for a vector $\nu$, $\nu_{(j)}$ denotes its $j$th element and $|\nu|$ denotes its Euclidean norm. We use $a_n \sim b_n$ to mean that there are constants $0 < A < B < \infty$ such that $A \leq a_n/b_n \leq B$ for all $n$.

CONDITION 2. The conditional distribution of $T$, given $X = x$, has a bounded density $f_{T|X} : 0 < c_1 \leq f_{T|X}(t|x) \leq c_2 < \infty$ uniformly in $x$ and $t$ for some positive constants $c_1$ and $c_2$.

CONDITION 3. $X_{(0)} = 1$, $E(X_{(j)}|T) = 0$ and $P\{|X_{(j)}| < M\} = 1$ for some $M < \infty$, $j = 1, \ldots, p$. There exist two positive definite matrices $\Sigma_1$ and $\Sigma_2$ such that $\Sigma_1 \leq \mathbf{Var}(X|T) \leq \Sigma_2$ uniformly in $T$, where $\mathbf{Var}(X|T)$ denotes the conditional covariance matrix of $X$ given $T$.

CONDITION 4. The $e_i$ are i.i.d. and have a density function $f_e$ that is continuous at 0 with $0 < f_e(0) < \infty$.

With an appropriate choice of $k_n$ to balance bias and variance, the proposed polynomial spline estimator and its derivatives attain the optimal convergence rates established by [31], as shown in the following theorem:

THEOREM 1. *Assume Conditions 1–4. Suppose $k_n \sim n^{1/(2r+1)}$ and $r > 1/2$. Then for $j = 0, 1, \ldots, p$,*

$$(2.3) \quad \frac{1}{n}\sum_{i=1}^{n}(\hat{\beta}_j^{(k)}(T_i) - \beta_j^{(k)}(T_i))^2 = O_p(n^{-2(r-k)/(2r+1)}), \qquad k = 0, 1, \ldots, m.$$

If piecewise linear splines are used and the $\beta_j(t)$ have bounded second order derivatives, then Theorem 1 gives the rate of convergence as $n^{-2/5}$.

REMARK 1. Considerable effort has been directed at relaxing the 'nonstochastic quasi-uniform knots assumption' in least squares regression. For example, Mao and Zhao [27] and Stone and Huang [33] investigated free-knot splines in which both the knot locations and spline coefficients were treated as unknown parameters. Extending the methodology and theory of free-knot splines to quantile regression is beyond the scope of this paper.

2.3. *Implementation.* When the B-spline basis is given, computations can be performed using standard quantile regression algorithms. As for selecting the order and knots for the splines, we propose a simple semi-automatic stepwise algorithm: users determine the order and specify a set of potential knots and the knots are chosen automatically from the potential knot set. The algorithm adds or deletes knots iteratively using Rao and Wald statistics to avoid additional model fitting. It allows different knots to be chosen for different coefficients. The algorithm terminates when no further addition or deletion occurs or when the model at the previous iteration is the same as the one at the current iteration. Similar algorithms are found in [26] and [32]. Upon termination, the algorithm finds the best fitting model using an adapted Schwarz-type Information Criterion (SIC),

$$SIC(M^i) = \log\left(\sum_{j=1}^{n} \rho_\tau(r_j^i)\right) + 0.5\log(n)p_n^i/n, \quad (2.4)$$

where $M^i$ denotes the $i$th interim model, $r_j^i$ is the $j$th residual from fitting $M^i$, $p_n^i$ is the number of variables in $M^i$ and $n$ is the sample size. For algorithmic convenience, we use a truncated power basis. We refer to Chapter 3.2 of [22] for more details.

2.3.1. *Order of splines to fit.* We suggest using lower order splines, such as linear ($m = 1$) and quadratic splines ($m = 2$), for practical reasons. Since the effect of the splines on the model is multiplicative, higher order splines would induce complicated interactions and collinearity among the variables in the model. For example, even the simplest cubic splines (one-piece cubic polynomials) would induce interactions of the form $x_j t$, $x_j t^2$ and $x_j t^3$. We prefer linear splines in particular since with them, the coefficient estimates are easier to interpret. Moreover, linear splines have an optimal property [24].

2.3.2. *Potential knot set.* Users provide either a potential knot set or the number of potential knots. In the latter case, the algorithm chooses as potential knots the user-provided number of equispaced knots. A helpful guideline of $\min(4n^{1/5}, n/4, N, 30)$ is found in [32] for the number of potential knots, where $N$ is the number of distinct data points.

**3. Model assessment.** Model assessment takes the form of a hypothesis test examining whether all of the coefficients are constant. Under Condition 1, this hypothesis test can be represented as

$$H_0 : q_\tau(t, x) \in \mathcal{L} \quad \text{vs.} \quad H_1 : q_\tau(t, x) \in \mathcal{G} \backslash \mathcal{L}, \quad (3.1)$$

where $\mathcal{L}$ and $\mathcal{G}$ are sets of functions defined as

$\mathcal{L} = \{q(t,x) \,|\, q(t,x) = b_0 x_0 + \cdots + b_p x_p,\ \text{for some constants } b_j,\ j = 0, \ldots, p\},$

$\mathcal{G} = \{q(t,x) \,|\, q(t,x) = \beta_0(t) x_0 + \cdots + \beta_p(t) x_p,\ \text{for } \beta_j(t) \in \mathcal{H}_r,\ j = 0, \ldots, p\}.$

We consider a transformation of $\Pi_{k_n}(t, x)$: $\Gamma_{k_n}(t, x) = \mathbf{A}_{k_n} \Pi_{k_n}(t, x) = (x^{\mathrm{T}}, \Gamma_{2,k_n}(t, x)^{\mathrm{T}})^{\mathrm{T}}$ for some transformation matrix $\mathbf{A}_{k_n}$. We denote the vector of coefficients of $\Gamma_{k_n}(t, x)$ by $\tilde{\boldsymbol{\xi}}_{k_n} = (\tilde{\boldsymbol{\xi}}_{1,k_n}^{\mathrm{T}}, \tilde{\boldsymbol{\xi}}_{2,k_n}^{\mathrm{T}})^{\mathrm{T}}$, where $\tilde{\boldsymbol{\xi}}_{1,k_n}$ and $\tilde{\boldsymbol{\xi}}_{2,k_n}$ correspond to $x$ and $\Gamma_{2,k_n}(t, x)$, respectively. Then the minimization problem is

$$\min_{\boldsymbol{\xi}_{k_n} \in \mathbb{R}^{p_{k_n}}} \sum_{i=1}^{n} \rho_\tau(Y_i - \boldsymbol{\xi}_{k_n}^{\mathrm{T}} \Gamma_{k_n}(T_i, X_i)) \quad (3.2)$$

and (3.1) can be represented as $H_0 : \tilde{\boldsymbol{\xi}}_{2,k_n} = 0$ versus $H_1 : \tilde{\boldsymbol{\xi}}_{2,k_n} \neq 0$. We denote the estimates of $\tilde{\boldsymbol{\xi}}_{k_n}$ obtained under $H_0$ and $H_1$ by $\hat{\boldsymbol{\xi}}_{k_n} = (\hat{\boldsymbol{\xi}}_1^{\mathrm{T}}, \mathbf{0}^{\mathrm{T}})^{\mathrm{T}}$ and $\bar{\boldsymbol{\xi}}_{k_n} = (\bar{\boldsymbol{\xi}}_{1,k_n}^{\mathrm{T}}, \bar{\boldsymbol{\xi}}_{2,k_n}^{\mathrm{T}})^{\mathrm{T}}$, respectively.

The proposed Rao-score-type test uses the score of (3.2) evaluated at the estimates of $\tilde{\boldsymbol{\xi}}_{k_n}$ under $H_0$. Let $s_{k_n} = (k_n/n)^{1/2} \sum_{i=1}^n \varphi_\tau(Y_i - \hat{\boldsymbol{\xi}}_1^{\mathrm{T}} X_i) \Gamma_{2,k_n}(T_i, X_i)$, where $\varphi_\tau(e) = \tau I(e > 0) + (\tau - 1) I(e < 0)$, the derivative of $\rho_\tau(e)$. Define $\sigma^2 = E(\varphi_\tau(e)^2)$ and $\mathbf{Q}_n = \sum_{i=1}^n \Gamma_{k_n}(T_i, X_i) \Gamma_{k_n}(T_i, X_i)^{\mathrm{T}}$. Let $\mathbf{Q}_{n(ij)}$, $i, j = 1, 2$, denote the $ij$th block of $\mathbf{Q}_n$ such that $\mathbf{Q}_{n(12)} = \sum_{i=1}^n X_i \Gamma_{2,k_n}(T_i, X_i)^{\mathrm{T}}$, for example. The test statistic is

$$(3.3) \qquad w_{k_n} = (1/\sigma)(k_n \mathbf{Q}_n^{(22)}/n)^{1/2} s_{k_n},$$

where $\mathbf{Q}_n^{(22)} = (\mathbf{Q}_{n(22)} - \mathbf{Q}_{n(21)} \mathbf{Q}_{n(11)} \mathbf{Q}_{n(12)})^{-1}$.

The test statistic has a simple asymptotic null distribution. Suppose that $k_n = k_0$, a positive constant. Under Conditions 1–4, for any $q_\tau(t, x) \in \mathcal{L}$, $|w_{k_0}|^2 \longrightarrow^D \chi^2_{(p+1)(k_0+1)}$ as $n \to \infty$ [23]. For increasing $k_n$, we have the following theorem:

THEOREM 2. *Assume Conditions 1–4. Suppose $\lim_{n \to \infty} k_n^2 n^{\delta-1} = 0$ for some $0 < \delta < 1$. Then for any $q_\tau(t, x) \in \mathcal{L}$,*

$$\frac{|w_{k_n}|^2 - (p_{k_n} - p - 1)}{\sqrt{2(p_{k_n} - p - 1)}} \longrightarrow^D N(0, 1) \qquad as\ p_{k_n} \to \infty.$$

If $k_n$ is bounded, then one uses the chi-square distribution as the limiting distribution. If $k_n$ is unbounded, then the standard normal distribution is used for the appropriately standardized test statistic.

The proposed test is simple and straightforward: it does not require estimation of $f_e(0)$ and the asymptotic null distribution of the test statistic is tractable. One might consider a Wald-type test as a natural choice for $\tilde{\boldsymbol{\xi}}_{2,k_n}$ under question. However, this would require estimation of $f_e(0)$ and would therefore be less desirable.

For purpose of comparison, we consider the following likelihood-type-test. We use the objective function in (3.2) in place of the likelihood function and define the test statistic as

$$l_{k_n} = 2 \left\{ \sum_{i=1}^n \rho_\tau(Y_i - \hat{\boldsymbol{\xi}}_1^{\mathrm{T}} X_i) - \sum_{i=1}^n \rho_\tau(Y_i - \bar{\boldsymbol{\xi}}_{1,k_n}^{\mathrm{T}} X_i - \bar{\boldsymbol{\xi}}_{2,k_n}^{\mathrm{T}} \Gamma_{2,k_n}(T_i, X_i)) \right\}.$$

As in the parametric case, the likelihood-type test also requires estimation of $f_e(0)$ [23]. Following [1, 2], we simulate the null distribution via the residual bootstrap.

3.1. *Heteroscedastic errors.* Heteroscedastic error models are of great interest in quantile regression. The Rao-score-type test is applicable to a scale family of linear heteroscedastic models such that

$$Y_i = \boldsymbol{\theta}^{\mathrm{T}} X_i + s(T_i, X_i) e_i, \tag{3.4}$$

where $s(T_i, X_i) > 0$ is a scale function that is consistently estimable and the $e_i$ are i.i.d.

A heuristic argument for this is as follows. The test only requires fitting the null model and the distribution of the test statistic relies only on the quantile estimator under the null. Therefore, any quantile estimator with the same first order representation under the null leads to the same limiting distribution for the test statistic. This implies that if we can obtain a $\sqrt{n}$-consistent estimator $\hat{\sigma}_i = \hat{s}(T_i, X_i)$ and if the weighted quantile estimator minimizing $\sum_{i=1}^{n} \hat{\sigma}_i^{-1} \rho_\tau(Y_i - \boldsymbol{\theta}^{\mathrm{T}} X_i)$ is first-order equivalent to the quantile estimator with $s(T_i, X_i)$, then we can apply the Rao-score-type test to $((T_i, X_i)/\hat{\sigma}_i, Y_i/\hat{\sigma}_i)$ as if the rescaled observations were i.i.d. A primary example of interest is $s(T_i, X_i) = (T_i, X_i^{\mathrm{T}})\boldsymbol{\gamma}$. We can obtain a $\sqrt{n}$-consistent estimator $\hat{\boldsymbol{\gamma}}$ from the absolute residuals of median regression and Koenker and Zhao [25] showed that the weighted quantile estimator minimizing $\sum_{i=1}^{n} \hat{\sigma}_i^{-1} \rho_\tau(Y_i - \boldsymbol{\theta}^{\mathrm{T}} X_i)$ with $\hat{\sigma}_i = (T_i, X_i^{\mathrm{T}})\hat{\boldsymbol{\gamma}}$ is first-order equivalent to the quantile estimator with $\boldsymbol{\gamma}$. For general $s(T_i, X_i)$, Zhao [36] showed that $s(T_i, X_i)$ can be consistently estimated using a nearest neighbor approach and that the resulting weighted quantile estimator is first-order equivalent. On the contrary, the likelihood-ratio-type test requires the alternative model to be fitted; hence, whether (or how) it can be applied to (3.4) is unclear.

3.2. *Test consistency.* Both tests encounter a problem of inconsistency when $k_n$ is fixed since the constancy of the coefficient functions is tested via the constancy of the approximating polynomial splines. To be specific, we consider a set of functions $\mathbb{L}_{k_n}$ such that

$$\mathbb{L}_{k_n} = \{q(t,x) \in \mathcal{G} \mid q(t,x) = \tilde{\boldsymbol{\xi}}_1^{\mathrm{T}} x - R_{g,k_n}(t,x)\},$$

where $\tilde{\boldsymbol{\xi}}_{k_n} = (\tilde{\boldsymbol{\xi}}_1^{\mathrm{T}}, \mathbf{0}^{\mathrm{T}})^{\mathrm{T}}$ and the $R_{g,k_n}(t,x)$ are given in Lemma A.1(ii). We note that $\mathbb{L}_{k_n}$ is a subset of $\mathcal{G}$ for which the approximations in the given spline space are functions of $x$ with constant coefficients. When $k_n = k_0$, testing (3.1) means testing $H_0: q_\tau(t,x) \in \mathbb{L}_{k_0}$ versus $H_1: q_\tau(t,x) \in \mathcal{G} \setminus \mathbb{L}_{k_0}$. Note that $\mathcal{L} \subset \mathbb{L}_{k_0}$, while $\mathbb{L}_{k_0} \setminus \mathcal{L} \neq \varnothing$. Also, any element of $\mathbb{L}_{k_0} \setminus \mathcal{L}$ is a null model and the proposed test will not be able to reject $H_0$, even as $n \to \infty$. The set $\mathbb{L}_{k_0} \setminus \mathcal{L}$ constitutes a class of nonconstant functions in $t$ associated with the fixed number of knots $k_0$, against which the proposed test is not consistent.

When $k_n$ is allowed to grow with $n$, we can consider a sequence $\{\mathbb{L}_{k_n}\}$ with $\mathbb{L}^* = \liminf_{n \to \infty} \mathbb{L}_{k_n}$. With respect to this sequence, (3.1) can be restated as $H_0: q_\tau(t, x) \in \mathbb{L}^*$ versus $H_1: q_\tau(t, x) \in \mathcal{G} \backslash \mathbb{L}^*$. The following lemma shows that the test is consistent for all alternatives when $k_n$ increases.

LEMMA 1. *Under Condition 1, $\mathbb{L}^* = \mathcal{L}$.*

Fortunately, these issues do not pose a real concern in practice as the knots are chosen adaptively. For the null, the number of adaptively chosen knots will not grow with $n$ and the results of [23] ensure that the proposed test is level-appropriate. For the alternative, the number of adaptively chosen knots will grow with $n$ and we have the above lemma for consistency of the test along with Theorem 2 to ensure that the test is level-appropriate.

## 4. Empirical results.

4.1. *Monte Carlo studies.* We conducted Monte Carlo studies for the following i.i.d. and heteroscedastic error models in order to evaluate the proposed test relative to the likelihood-ratio-type test:

$$\text{M1:} \quad y_i = \gamma_1 + b(t_i)x_i + (e_i - F_e^{-1}(\tau)),$$
$$\text{M2:} \quad y_i = \gamma_2 + bx_i + (\gamma_3 t_i + \gamma_4 x_i)(e_i - F_e^{-1}(0.5)),$$

where $F_e(\cdot)$ denotes the distribution function of the $e_i$, the $x_i$ are truncated standard normal variables, that is ($x_i = sign(z_i) \min(c, |z_i|)$ with $z_i$ standard normal and $c > 0$), the $\gamma_j$ are constants, $t_i \sim U(0.5, 1.5)$ for all $i$ and $t_i$ is independent of $x_i$ and $e_i$ for all $i$. Three different distributions were used for the $e_i$: a standard normal, a $\chi^2$ with 1 degree of freedom and a $t$ with 3 degrees of freedom. M1 represents a homoscedastic error model that is a null model when $b(t_i) = b$ and an alternative model when $b(t_i)$ is a nonconstant function of $t$. Four different alternatives were considered for $b(t_i)$: a linear, quadratic, sine and log function of $t$. M2 represents a heteroscedastic error model that is a null model at the conditional median ($\tau = 0.5$).

The Monte Carlo studies were based on 500 data sets, each a sample of size 200. We used piecewise linear splines with an adaptively chosen number of uniform knots to fit a model to each simulated sample. The number of knots was chosen by the adapted SIC in (2.4). This closely resembles standard practice, in which the order of splines and knot placements are chosen nonadaptively, but the number of knots is chosen adaptively. Piecewise linear splines were used because they provided reasonable fits for the various alternatives considered here. For the likelihood-ratio-type test, 200 bootstrapped samples were used to simulate the null distribution of the test statistic. The performance of the tests was measured by power at the significance level $\alpha = 0.05$. As the knots were chosen adaptively, powers under the null, that is the type I error rates, were expected to be higher than the nominal significance level.

TABLE 1
*Power at $\alpha = 0.05$ when errors are i.i.d. (M1). The parenthetical numbers in the second column indicate the densities of the respective error distributions at the $\tau$th quantile. "RS" and "LR" indicate the Rao-score- and likelihood-ratio-type tests*

| $\tau$ | Error distribution | constant | | linear | | quadratic | | sine | | log | |
|---|---|---|---|---|---|---|---|---|---|---|---|
| | | RS | LR | RS | LR | RS | LR | RS | LR | RS | LR |
| 0.5 | $Z$ (0.40) | 0.062 | 0.072 | 0.71 | 0.80 | 0.69 | 0.80 | 0.84 | 0.90 | 0.75 | 0.86 |
| | $\frac{1}{\sqrt{2}}\chi_1^2$ (0.67) | 0.040 | 0.062 | 1 | 1 | 1 | 1 | 1 | 1 | 1 | 1 |
| | $\frac{1}{\sqrt{3}}t_3$ (0.64) | 0.054 | 0.068 | 0.96 | 0.98 | 0.94 | 0.99 | 0.99 | 1 | 0.98 | 0.99 |
| 0.9 | $Z$ (0.18) | 0.070 | 0.130 | 0.39 | 0.64 | 0.4 | 0.64 | 0.43 | 0.75 | 0.46 | 0.65 |
| | $\frac{1}{\sqrt{2}}\chi_1^2$ (0.09) | 0.082 | 0.164 | 0.16 | 0.33 | 0.16 | 0.31 | 0.19 | 0.41 | 0.17 | 0.32 |
| | $\frac{1}{\sqrt{3}}t_3$ (0.18) | 0.056 | 0.130 | 0.40 | 0.62 | 0.49 | 0.70 | 0.50 | 0.78 | 0.55 | 0.72 |

Table 1 summarizes the results for the homoscedastic error model. The likelihood-ratio-type test was slightly more powerful, but less likely to attain the nominal significance level, particularly when $\tau = 0.9$. The Rao-score-type test performed reliably with comparable power. Inherent to quantile regression, the performances of both tests depend on the error density at the quantile under consideration; power rises with increasing error density.

Table 2 summarizes the results for the heteroscedastic error median null model. The Rao-score-type test attains type I error rates close to the nominal significance level. The likelihood-ratio-type test was administered as if the errors were i.i.d. since it was not clear how to simulate the null distribution of the test statistic via the bootstrap. Quite predictably, it did not perform well.

4.2. *Example.* We illustrate the prescribed methodology using the forced expiratory volume (FEV) data from [29]. Measurements of FEV in liters, age ($T$) in years, height ($H$) in inches, sex ($S = 1$ for boys/$S = 0$ for girls) and smoking status ($SM = 1$ for a current smoker/$SM = 0$ otherwise) were collected for 654 children aged 3–19 who participated in the Childhood Respiratory Disease Study.

Of interest are the effect of cigarette smoking on FEV and the lower conditional quantiles of FEV as a gauge for poor pulmonary function. The proposed methodology enables us to investigate both research questions while allowing the coefficients to vary with age to reflect the possibility of age-dependent covariate effects and the lower quantiles to be related to the covariates in ways different from those for the mean. We modeled the first quartile ($q_{0.25}$) and the median ($q_{0.5}$) as follows:

(4.1) $\quad q_\tau = \beta_0(t) + \beta_1(t) \times S + \beta_2(t) \times H + \beta_3(t) \times (H \times S) + \beta_4(t) \times SM.$

We tried both piecewise linear and quadratic splines with a set of eleven equispaced knots, using the stepwise knot selection algorithm described in

Table 2
*Type I Error Rate at $\alpha = 0.05$ when errors are heteroscedastic (M2 with $\tau = 0.5$)*

| \multicolumn{2}{c}{$Z$} | \multicolumn{2}{c}{$\frac{1}{\sqrt{2}}\chi_1^2$} | \multicolumn{2}{c}{$\frac{1}{\sqrt{3}}t_3$} |
|---|---|---|---|---|---|
| RS | LR | RS | LR | RS | LR |
| 0.076 | 0.378 | 0.084 | 0.432 | 0.056 | 0.324 |

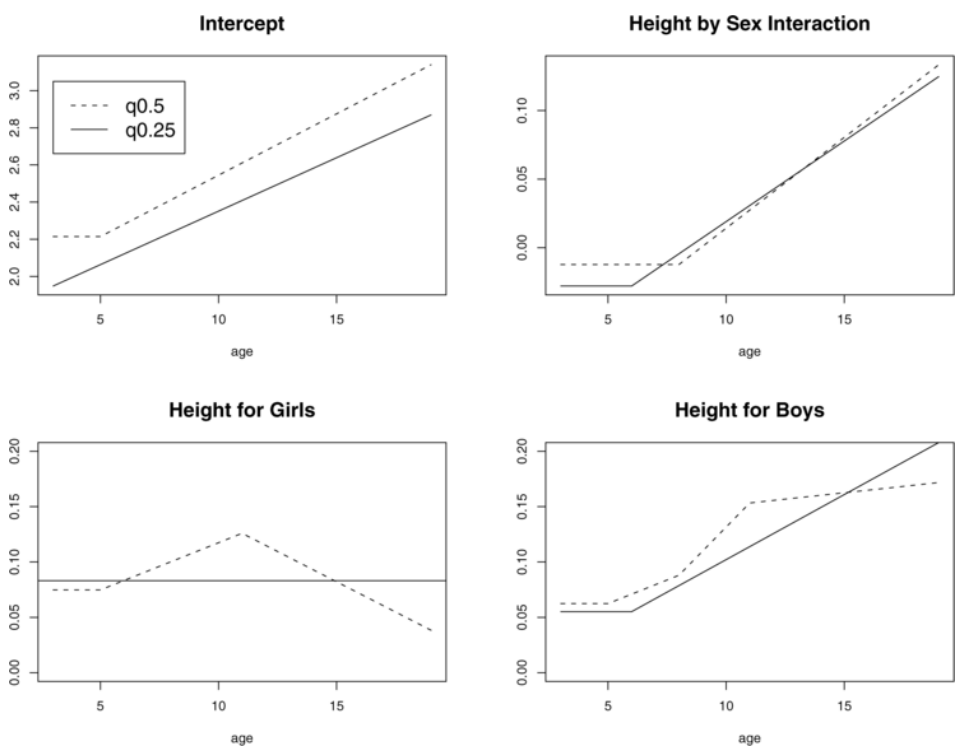

Fig. 1. *Coefficient estimates by piecewise linear splines.*

Section 2.3. Both splines provided similar fits, while the fit with piecewise linear splines was preferred by the adapted SIC. Hence, we shall discuss the results based on piecewise linear splines.

The estimates of $\beta_4(t)$ were small negative constants at both quantiles. This implies that the effect of smoking on FEV is not age dependent at these quantiles. This is not surprising since the current smokers in the data set were close to each other in age.

$\beta_0(t)$, $\beta_2(t)$ and $\beta_3(t)$ were estimated as varying with age (see Figure 1). At both quantiles, the estimated intercept increases linearly, implying that FEV increases linearly with age. The height coefficient for girls was esti-mated as a constant for $q_{0.25}$, while knots were chosen at 5 and 11 for $q_{0.5}$. For the coefficient of height by sex interaction, knots were chosen at 6 for $q_{0.25}$ and at 8 for $q_{0.5}$. The height coefficient estimates for boys in the bottom right panel were obtained as $\hat{\beta}_2(t) + \hat{\beta}_3(t)$.

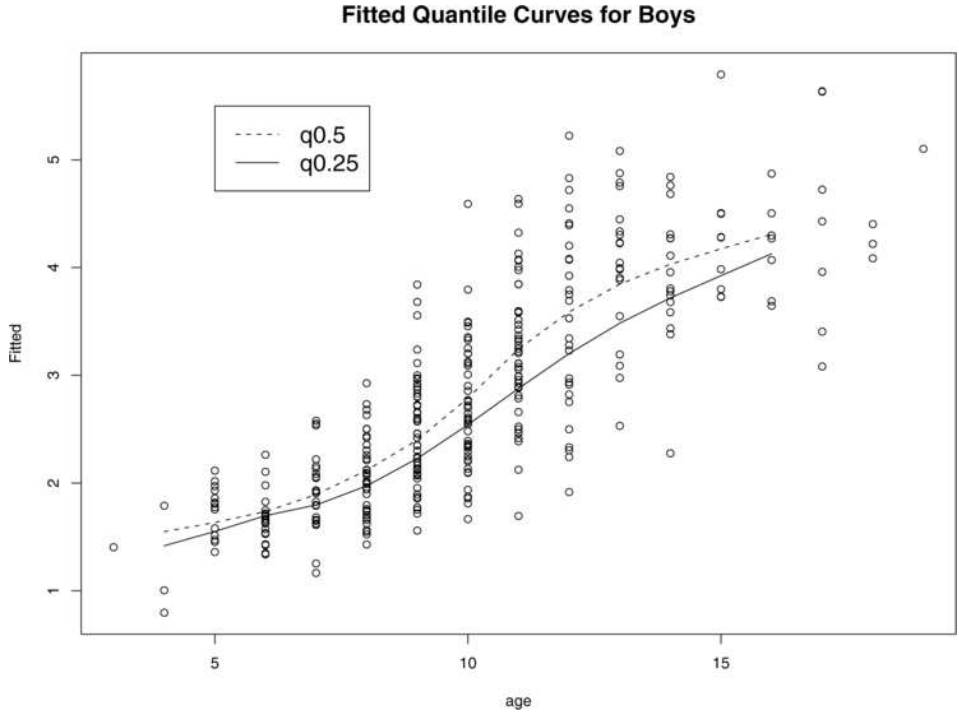

Fig. 2. *The fitted values are the quantile estimates obtained by plugging in the smoothed sample height average of each age group with the coefficient estimate. The curves are plotted over the credible age range where 90% of data are present.*

We focus on the estimates for ages from 5 to 15 in Figure 1 since little data is available outside the age range. The nonparallel intercept estimates and nonoverlapping coefficient estimates imply a different regression relationship at the first quartile than at the median. More specifically, the estimated intercept does not grow as rapidly at the first quartile as at the median. A similar pattern is shown in the height coefficient estimates for boys. This suggests that in boys, FEV does not increase with age as rapidly in the lower tail of its conditional distribution as in the middle, implying that in boys, the difference between height-adjusted typical and poor pulmonary function increases with age. It manifests through increasing variability in FEV over that age range, while variability in height does not necessarily increase (see Figure 2). This structural heterogeneity would not have been revealed by an ordinary regression.

We also tested for age-dependent effects of the covariates using the proposed Rao score test. We used linear splines with the same knots selected for the regressions. We rejected the constancy hypothesis with $p$-values near zero for both quantiles.

While FEV is known to increase with age, it is less clear which other coefficients change with age. In particular, do the coefficients of height and height by sex interaction change? To address this question, we detrended the data by subtracting age-and-gender-specific sample quartiles from the observations of respective age-gender groups. We then modeled the conditional quartiles of the detrended data as a function of height and height by sex interaction without an intercept. We used piecewise splines with the same knots selected for the regressions. We rejected the constancy hypothesis with a $p$-value $< 0.006$ for $q_{0.25}$ and a $p$-value near zero for $q_{0.5}$. The results agreed with the nonconstant coefficient estimates obtained from the regressions.

## APPENDIX

For simplicity, proofs are provided for uniform knots and $p = 2$, the subscript $k_n$ being suppressed if possible. The same proofs will work for quasi-uniform knots and/or $p > 2$. We adopt the following notation:

$$\mathbf{G}_n = (\pi(T_i), \ldots, \pi(T_n))(\pi(T_i), \ldots, \pi(T_n))^{\mathrm{T}},$$

$$\widetilde{\mathbf{X}}_n^{\mathrm{T}} = (\Pi(T_i, X_i), \ldots, \Pi(T_n, X_n))_{p_{k_n} \times n}, \qquad \mathbf{H}_n = \widetilde{\mathbf{X}}_n^{\mathrm{T}} \widetilde{\mathbf{X}}_n,$$

$$\mathbf{G} = E(\pi(T)\pi(T)^{\mathrm{T}}), \qquad \mathbf{H} = E(\Pi(T,X)\Pi(T,X)^{\mathrm{T}}), \qquad r_{ij} = r_{j,k_n}(T_i),$$

$$R_i = R_{g,k_n}(T_i, X_i), \qquad \hat{\boldsymbol{\theta}}^* = \mathbf{H}_n^{1/2} \hat{\boldsymbol{\theta}}, \qquad \tilde{\boldsymbol{\theta}}^* = \mathbf{H}_n^{1/2} \tilde{\boldsymbol{\theta}},$$

$$Z_i = \mathbf{H}_n^{-1/2} \Pi(T_i, X_i),$$

where $\mathbf{H}_n^{-1/2}$ denotes the Moore inverse of $\mathbf{H}_n^{1/2}$. Also, for a matrix $\mathbf{D}$, let $\mathbf{D}_{(i,j)}$ denote its element at the $i$th row and $j$th column. For example, $\mathbf{H}_{n(aN+c, bN+d)} = \sum_{i=1}^n B_c(T_i) B_d(T_i) X_{i(a)} X_{i(b)}$ and $\mathbf{H}_{(aN+c, bN+d)} = E(B_c(T) B_d(T) X_{(a)} X_{(b)})$ for $a, b = 0, 1, 2$ and $1 \leq c, d \leq k_n + m + 1$.

We first present lemmas that are necessary to prove the theorems. Lemma A.1 follows directly from Corollary 6.21 of [30] and the proofs of Lemmas A.2–A.4 follow lines of argument similar to those used in the proofs of lemmas in [5]. Lemma A.5 follows directly from Lemma A.4 and Condition 3. Therefore, we omit their proofs.

LEMMA A.1. *Assume Conditions 1 and 3. For some constants $W_1$ and $W_1^*$ that depend on only (but not all of) $m$, $W_0$, $p$ and $M$:*

(i) there exist $\tilde{\theta}_j^{\mathrm{T}} \pi_{k_n}(t)$, $j = 0, \ldots, p$, such that $\beta_j(t) = \tilde{\theta}_j^{\mathrm{T}} \pi_{k_n}(t) - r_{j,k_n}(t)$ and $\sup_{t \in [0,1]} |r_{j,k_n}(t)| \le W_1 k_n^{-r}$;

(ii) there exists $\tilde{\boldsymbol{\theta}}_{k_n}^{\mathrm{T}} \Pi_{k_n}(t,x)$ such that $q_\tau(t,x) = \tilde{\boldsymbol{\theta}}_{k_n}^{\mathrm{T}} \Pi_{k_n}(t,x) - R_{g,k_n}(t,x)$ and $\sup_{(t,x) \in [0,1] \times R^{p+1}} |R_{g,k_n}(t,x)| \le W_1^* k_n^{-r}$.

LEMMA A.2. *Assume Conditions 1–3. For all $k_n$, the eigenvalues of $k_n \mathbf{G}$ and $k_n \mathbf{H}$ are bounded and depend only on $m$.*

LEMMA A.3. *Assume Conditions 1–3. Suppose $\lim_{n \to \infty} k_n n^{\delta - 1} = 0$ for some $0 < \delta < 1$. Except on an event which depends on $\{(T_1, X_1), (T_2, X_2), \ldots, (T_n, X_n)\}$ and whose probability tends to zero with increasing $n$, $|\mathbf{H}_{n(aN+c,bN+d)}/n - \mathbf{H}_{(aN+c,bN+d)}| = O_p(c_n n^{-1/2} k_n^{1/2})$, uniformly in $a, b, c$ and $d$, where $\{c_n\}$ is a nondecreasing sequence of positive numbers for which $c_n^2/n k_n^{-1} \to 0$ and $c_n^2/\log k_n \to \infty$.*

LEMMA A.4. *Assume Conditions 1–3. Suppose $\lim_{n \to \infty} k_n n^{\delta - 1} = 0$ for some $0 < \delta < 1$. Then the eigenvalues of $k_n \mathbf{G}_n/n$ and $k_n \mathbf{H}_n/n$ are bounded in probability.*

LEMMA A.5. *Assume Conditions 1–3. Suppose $\lim_{n \to \infty} k_n n^{\delta - 1} = 0$ for some $0 < \delta < 1$. Then with probability one, $\max_{1 \le i \le n} |Z_i| = O(\sqrt{k_n/n})$.*

LEMMA A.6. *Under the conditions of Theorem 1, $(\hat{\boldsymbol{\theta}} - \tilde{\boldsymbol{\theta}})^{\mathrm{T}} \mathbf{H}_n (\hat{\boldsymbol{\theta}} - \tilde{\boldsymbol{\theta}}) = O_p(k_n)$.*

PROOF. Since $(\hat{\boldsymbol{\theta}} - \tilde{\boldsymbol{\theta}})^{\mathrm{T}} \mathbf{H}_n (\hat{\boldsymbol{\theta}} - \tilde{\boldsymbol{\theta}}) = (\hat{\boldsymbol{\theta}}^* - \tilde{\boldsymbol{\theta}}^*)^{\mathrm{T}} (\hat{\boldsymbol{\theta}}^* - \tilde{\boldsymbol{\theta}}^*)$, it suffices to show that $|\hat{\boldsymbol{\theta}}^* - \tilde{\boldsymbol{\theta}}^*|^2 = O_p(k_n)$. As $\hat{\boldsymbol{\theta}}^* = \arg\min_{\boldsymbol{\theta} \in \mathbb{R}^{p_{k_n}}} \sum_{i=1}^n \rho_\tau(e_i - (\boldsymbol{\theta} - \tilde{\boldsymbol{\theta}}^*)^{\mathrm{T}} Z_i - R_i)$ from (2.2), with probability 1 we have

$$(\mathrm{A.1}) \quad \sum_{i=1}^n \rho_\tau(e_i - (\hat{\boldsymbol{\theta}}^* - \tilde{\boldsymbol{\theta}}^*)^{\mathrm{T}} Z_i - R_i) = \inf_{\boldsymbol{\theta} \in \mathbb{R}^{p_{k_n}}} \sum_{i=1}^n \rho_\tau(e_i - \boldsymbol{\theta}^{\mathrm{T}} Z_i - R_i).$$

Let $\Theta_L = \{\boldsymbol{\theta} | \boldsymbol{\theta} \in \mathbb{R}^{p_{k_n}}, |\boldsymbol{\theta}| \le L k_n^{1/2}\}$ for some constant $L > 0$. Then the proof is reduced to showing that a solution to the optimization problem in (A.1) lies in $\Theta_L$ for a sufficiently large $L$. The following lemma and the convexity of the minimization problem (A.1) complete the proof. $\square$

LEMMA A.7. *Assume the conditions of Theorem 1.*

(i) *For any sequence $\{L_n\}$ satisfying $1 \le L_n \le k_n^{\eta_0/10}$ for some $0 < \eta_0 < (r - 1/2)/(2r + 1)$, we have*

$$\sup_{|\boldsymbol{\theta}| \le L_n k_n^{1/2}} k_n^{-1} \left| \sum_{i=1}^n [\rho_\tau(e_i - Z_i^{\mathrm{T}} \boldsymbol{\theta} - R_i) - \rho_\tau(e_i - R_i) + Z_i^{\mathrm{T}} \boldsymbol{\theta}(\tau - I(e_i < 0)) \right.$$

$$- E(\rho_\tau(e_i - Z_i^{\mathrm{T}}\boldsymbol{\theta} - R_i) - \rho_\tau(e_i - R_i))]\bigg| = o_p(1).$$

(ii) *For any $\varepsilon > 0$, there exists $L := L_\varepsilon$ (sufficiently large) such that as $n \to \infty$,*

$$P\bigg\{k_n^{-1}\bigg(\inf_{|\boldsymbol{\theta}|=Lk_n^{1/2}}\sum_{i=1}^n[E(\rho_\tau(e_i - Z_i^{\mathrm{T}}\boldsymbol{\theta} - R_i) - \rho_\tau(e_i - R_i))]$$
$$- \bigg|\sum_{i=1}^n Z_i(\tau - I(e_i < 0))\bigg|\bigg) > 1\bigg\} > 1 - \varepsilon.$$

The proof of Lemma A.7 uses the results of Lemmas A.4–A.5 and a partition of the parameter space $\Theta_L$. It applies Bernstein's inequality to bound the probabilities over the partition. The method is similar to that used in Lemmas 3.2 and 3.3 of [14] and so we omit the details.

By Lemma A.7, there exists some $L := L_\varepsilon$ for any $\varepsilon > 0$ such that as $n \to \infty$, we have

$$P\bigg\{\inf_{|\boldsymbol{\theta}|=Lk_n^{1/2}}\sum_{i=1}^n \rho_\tau(e_i - Z_i^{\mathrm{T}}\boldsymbol{\theta} - R_i) > \sum_{i=1}^n \rho_\tau(e_i - R_i)\bigg\} > 1 - \varepsilon.$$

By Corollary 25 of [8], we have that as $n \to \infty$,

(A.2) $$P\bigg\{\inf_{|\boldsymbol{\theta}|\geq Lk_n^{1/2}}\sum_{i=1}^n \rho_\tau(e_i - Z_i^{\mathrm{T}}\boldsymbol{\theta} - R_i) > \sum_{i=1}^n \rho_\tau(e_i - R_i)\bigg\} > 1 - \varepsilon.$$

By (A.2), (A.1) implies that $P\{|\hat{\boldsymbol{\theta}}^* - \tilde{\boldsymbol{\theta}}^*| \leq Lk_n^{1/2}\} > 1 - \varepsilon$.

PROOF OF THEOREM 1. First, it follows from Lemma A.1(i) that

$$\frac{1}{n}\sum_{i=1}^n(\hat{\beta}_j(T_i) - \beta_j(T_i))^2 \leq \frac{2}{n}\sum_{i=1}^n(\pi(T_i)^{\mathrm{T}}\hat{\theta}_j - \pi(T_i)^{\mathrm{T}}\tilde{\theta}_j)^2 + \frac{2}{n}\sum_{i=1}^n r_{ij}^2$$
$$\leq \frac{2}{n}(\hat{\theta}_j - \tilde{\theta}_j)^{\mathrm{T}}\mathbf{G}_n(\hat{\theta}_j - \tilde{\theta}_j) + 2W_2^2 k_n^{-2r}.$$

As $k_n \sim n^{1/(2r+1)}$, by Lemma A.4, it suffices to show that $|\hat{\theta}_j - \tilde{\theta}_j|^2 = O_p(k_n^2/n)$. By Lemmas A.4 and A.6, $|\hat{\boldsymbol{\theta}} - \tilde{\boldsymbol{\theta}}|^2 = O_p(k_n^2/n)$, which implies that $|\hat{\theta}_j - \tilde{\theta}_j|^2 = O_p(k_n^2/n)$. □

PROOF OF LEMMA 1. As $\mathcal{L} \subset \mathbb{L}^*$ by definition, it suffices to show that $\mathbb{L}^* \subset \mathcal{L}$. Suppose that $\mathbb{L}^*/\mathcal{L} \neq \varnothing$. Let $q^*(t,x)$ denote an element of $\mathbb{L}^*/\mathcal{L}$ and let $\{k_{n*}\}$ denote a subsequence of $\{k_n\}$ such that $q^*(t,x) \in \mathbb{L}_{k_{n*}}$ for all $k_{n*}$. Consider $R_{q^*,k_{n*}}(t,x)$, where $R_{q^*,k_{n*}}(t,x)$ is given in Lemma A.1(ii) with respect to $q^*(t,x)$. We note that $q^* \in \mathbb{L}^*/\mathcal{L}$ implies $\lim_{k_{n*} \to \infty} R_{g,k_{n*}}(t,x) \neq 0$. This contradicts Lemma A.1(ii). $\square$

PROOF OF THEOREM 2. Let $\zeta_i = (\mathbf{Q}_n^{(22)})^{1/2}\Gamma_2(T_i, X_i)$. Define $a_{k_n} = 1/\sigma \sum_{i=1}^n \zeta_i \varphi_\tau(e_i)$ and $b_{k_n} = 1/\sigma \sum_{i=1}^n \zeta_i \{\varphi_\tau(e_i - (\hat{\boldsymbol{\xi}}_1 - \tilde{\boldsymbol{\xi}}_1)^T x_i) - \varphi_\tau(e_i))\}$. Then $w_{k_n} = a_{k_n} + b_{k_n}$. We have $\frac{|a_{k_n}|^2 - (p_{k_n} - p - 1)}{\sqrt{2(p_{k_n} - p - 1)}} \longrightarrow^D N(0,1)$ from Theorem 4.1 of [28] by replacing the continuity condition on $\varphi_\tau(u)$ in the theorem with Condition 2. As $|a_{k_n}^T b_{k_n}| \leq |a_{k_n}||b_{k_n}|$ and $p_{k_n} - p - 1 = (p+1)(k_n + m)$, it suffices to show that $|b_{k_n}| = o_p(\sqrt{k_n})$. The following lemma completes the proof. $\square$

LEMMA A.8. *Assuming the conditions of Theorem 2, $|b_{k_n}| = o_p(\sqrt{k_n})$.*

PROOF. First, $F_e(t) - F_e(0) = f_e(0)t + o(|t|)$ and $E(\varphi_\tau(s+t) - \varphi_\tau(s))^2 = f_e(0)t + o(|t|)$ for sufficiently small $|t|$. Also, from Lemma A.5, $\max_i |\zeta_i| = O(\sqrt{k_n/n})$ with probability one.

Consider a function $b(\eta) = 1/\sigma \sum_{i=1}^n \zeta_i \{\varphi_\tau(e_i - x_i^T \eta/\sqrt{n}) - \varphi_\tau(e_i)\}$ for some $\eta \in \mathbb{R}^{p+1}$. For $b_{(j)}(\eta)$, $j = 1, \ldots, (p+1)(k_n + m)$, we have

$$Var(b_{(j)}(\eta)) \leq 1/\sigma^2 \sum_i^n E\{\varphi_\tau(e_i - x_i^T \eta/\sqrt{n}) - \varphi_\tau(e_i)\}^2 \zeta_{i(j)}^2 = O(k_n/\sqrt{n})$$

and $1/\sigma \max_i |\zeta_{i(j)}\{\varphi_\tau(e_i - x_i^T\eta/\sqrt{n}) - \varphi_\tau(e_i)\}| = O(\sqrt{k_n/n})$. It follows from Bernstein's inequality that for all $L > 0$, $\varepsilon > 0$ and any $|\eta| < L$,

$$P\{|b_{(j)}(\eta) - E(b_{(j)}(\eta))| > \varepsilon\} \leq 2\exp(-\sqrt{n}\varepsilon^2/(2k_n + 2\sqrt{k_n}/3)).$$

As the number of $b_{(j)}(\eta)$'s to consider is $(p+1)(k_n + m)$ and by the arbitrariness of $L > 0$ and $\varepsilon > 0$, we have

$$P\Big\{\sup_{|\eta| \leq L} |b(\delta) - E(b(\eta))| > \varepsilon\Big\}$$

$$\leq 2\exp(-\sqrt{n}\varepsilon^2/(2k_n + 2\sqrt{k_n}/3) + \log((p+1)(k_n + m))).$$

Under the conditions of Theorem 2, $nk_n^{-2}/(\log k_n)^2 \to \infty$ and the right-hand side converges to zero as $n \to \infty$. On the other hand, $|E(b(\eta))| = O(\sqrt{k_n/n})$, except for an event which depends on $\{(T_1, X_1), (T_2, X_2), \ldots, (T_n, X_n)\}$ and whose probability tends to zero with increasing $n$. Since $|\hat{\boldsymbol{\theta}}_{1k_n}| = O_p(1/\sqrt{n})$, the lemma follows. $\square$

CENTER FOR EPIDEMIOLOGY AND BIOSTATISTICS
CINCINNATI CHILDREN'S HOSPITAL MEDICAL CENTER, MLC 5041
3333 BURNED AVE.
CINCINNATI, OHIO 45229-3039
USA
E-MAIL: miok.kim@cchmc.org